\documentclass[12pt]{article}
\usepackage{amssymb,latexsym}
\usepackage{graphicx,psfrag}
\usepackage{amsmath}
\newcommand{\qed}{\nobreak \ifvmode \relax \else
\ifdim\lastskip<1.5em \hskip-\lastskip
\hskip1.5em plus0em minus0.5em \fi \nobreak
\vrule height0.75em width0.5em depth0.25em\fi}

\makeatletter

\renewcommand\section{\@startsection {section}{1}{\z@}
{-30pt \@plus -1ex \@minus -.2ex}
{2.3ex \@plus.2ex}
{\normalfont\normalsize\bfseries}}

\renewcommand\subsection{\@startsection{subsection}{2}{\z@}
{-3.25ex\@plus -1ex \@minus -.2ex}
{1.5ex \@plus .2ex}
{\normalfont\normalsize\bfseries}}

\renewcommand{\@seccntformat}[1]{\csname the#1\endcsname. }

\makeatother

\begin{document}

\title{On large primitive subsets of $\{1,2,\ldots,2n\}$}
\author{Sujith Vijay}
\date{}
\maketitle

\vskip 10pt

\centerline{\bf Abstract}

\vskip 10pt

\noindent

A subset of $\{1,2,\ldots,2n\}$ is said to be {\em{primitive}} if it does not contain any pair of elements $(u,v)$ such that $u$ is a divisor of $v$. Let $D(n)$ denote the number of primitive subsets of $\{1,2,\ldots,2n\}$ with $n$ elements. Numerical evidence suggests that $D(n)$ is roughly $(1.32)^n$. We show that for sufficiently large $n$, $$(1.303...)^n < D(n) < (1.408...)^n$$ 

\section*{\normalsize 1. Introduction}

A well-known application of the pigeonhole principle is the fact that any subset of  $\{1,2,\ldots,2n\}$ with more than $n$ elements must contain a pair of integers such that one divides the other. This is also famous in mathematical folklore as a ``recruiting problem" that Paul Erd\H{o}s asked young Lajos Po\'{s}a (see \cite{Aig98}). The trick is to assign each number in the subset to the pigeonhole corresponding to its largest odd divisor, whereupon the proof is immediate. It is also clear that the result is sharp, since $\{n+1,n+2,\ldots,2n\}$ is a  {\em{primitive}} (i.e., divisor-free) subset with $n$ elements. Various generalizations of this problem (see \cite{Heg06}, \cite{Leb77} and \cite{Vij06}) have also been studied. \\

Let $D(n)$ denote the number of primitive subsets of $\{1,2,\ldots,2n\}$ with $n$ elements. The first ten values of $D(n)$ are $2, 2, 3, 5, 4, 6, 12, 10, 14$ and $26$. The sequence is listed in the Online Encyclopedia of Integer Sequences (OEIS) as A174094. A natural question now arises. How fast does $D(n)$ grow? Numerical evidence suggests that $D(n)$ is roughly $(1.32)^n$. The purpose of this article is to derive upper and lower bounds on $D(n)$. \\

We note that the related function $D^{*}(n)$, denoting the number of all primitive subsets of $\{1,2,\ldots,2n\}$ (i.e, not necessarily with $n$ elements) has been studied by Erd\H{o}s and Cameron (see \cite{Cam90}). They have shown that $$(2.43...)^n < D^{*}(n) < (2.55...)^n$$ 

Recently, Angelo \cite{Ang18} established the existence of $$\lim_{n \rightarrow \infty} (D^{*}(n))^{1/n}$$

\section*{\normalsize 2. Upper Bound for $D(n)$}

Consider the partial order on the first $2n$ positive integers defined as follows: $a \preceq b$ if and only if $b/a$ is a power of $2$. Note that the first $2n$ positive integers can be partitioned into $n$ chains under this partial order. The least elements of the chains are precisely the odd integers $1,3,\ldots,2n-1$. \\

For brevity, we will refer to $n$-element primitive subsets as large primitive subsets (LPS). Clearly, every LPS of $\{1,2,\ldots,2n\}$ must contain exactly one element from each of these chains.  Thus an obvious upper bound is the product of the cardinalities of these chains. Since there are $\lfloor n/2 \rfloor$ chains of size $1$,  $\lfloor n/4 \rfloor$ chains of size $2$, $\lfloor n/8 \rfloor$ chains of size $3$ and so on, we have 
$$D(n) < 2^{\lfloor n/4 \rfloor} 3^{\lfloor n/8 \rfloor} \cdots < 2^{n/4} 3^{n/8} \cdots$$
$$\frac{\log_2(D(n))}{n} < \sum_{k=2}^{\infty} \frac{\log_2(k)}{2^k} = 0.7326...$$
$$D(n) < (1.661...)^n$$

We will improve this bound to show that $D(n) < (1.408...)^n$. Our approach will be to show that many integers either belong to all LPS or belong to none of them. We will color an integer green if it belongs to every LPS, red if it belongs to no LPS, and blue otherwise. Note that every LPS must contain all odd integers in $[n+1,2n]$, since they are the only elements in their chain. Thus all odd integers in $[n+1,2n]$ are green. \\

\noindent {\bf {Lemma 1}} Let $q \le 2n/3$ be an odd integer. Then $[n+1,2n]$ contains an odd multiple of $q$. \\

\noindent {\em {Proof.}} If $q \le  n/2$, the interval $[n+1,2n]$ has length at least $2q$, so it must contain two consecutive multiples of $q$, one of which must be odd. If $n/2 < q \le 2n/3$, the odd integer $3q$ belongs to $[n+1,2n]$. $\hfill \square$ \\

It follows from Lemma 1 that no LPS can contain any odd integer less than $2n/3$. Thus all odd integers in $[1,2n/3]$ are red. In particular, all odd integers in $(n/2,2n/3]$ are red. Let $q$ be such an integer. Since $4q > 2n$, the only elements in the chain containing $q$ are $q$ and $2q$. Since $q$ belongs to no LPS, $2q$ must belong to every LPS. Thus all integers congruent to $2$ modulo $4$ in $[n+1,4n/3]$ are green. \\

\noindent {\bf {Lemma 2}} Let $q$ be an odd integer belonging to any of the intervals $I_1=[1,2n/21], I_2=(n/10,2n/15]$ or $I_3=(n/6,2n/9]$. Then $[n+1,4n/3]$ contains an odd multiple of $2q$. \\

\noindent {\em {Proof.}} If $q \in I_2$, we have $10q \in [n+1,4n/3]$. Similarly, if $q \in I_3$, we have $6q \in [n+1,4n/3]$. \\

Let $q \in I_1$ be an odd integer. If $q \le  n/12$, the interval $[n+1,4n/3]$ has length at least $4q$, so it must contain two consecutive multiples of $2q$, one of which must be an odd multiple. If $n/12 < q \le 2n/21$, we have $14q \in [n+1,4n/3]$. $\hfill \square$ \\

It follows from Lemma 2 that all integers congruent to $2$ modulo $4$ in $[1,4n/21], (n/5, 4n/15]$ and $(n/3, 4n/9]$ are red. All the remaining elements are colored blue. (Some of them can possibly be colored red or green by more sophisticated arguments, but we do not advance them.) \\

We now count the number of blue elements in the chain containing $q$ for each odd integer $q \in [1,2n]$. As we have already seen, there are no blue elements for $q > n$ and $q \in (n/2,2n/3]$. \\

Let $J_1 = (2n/3,n], J_2 = (n/4,n/2], J_3 = (2n/9,n/4], J_4 = (n/6,2n/9], J_5 = (2n/15,n/6], J_6 = (n/8,2n/15], J_7 = (n/10,n/8], J_8 = (2n/21,n/10]$ and $J_9 = (n/16,2n/21]$. It follows from our discussions that there are two blue elements in the chain containing $q$ for $q \in J_1 \cup J_2 \cup J_4 \cup J_6$. Similarly, there are three blue elements in the chain containing $q$ for $q \in J_3 \cup J_5 \cup J_7 \cup J_9$, and four blue elements in the chain containing $q$ for $q \in J_8$. Additionally, there are $k$ blue elements in each chain for $q \in (n/2^{k+1},n/2^k]$ for  each $k \ge 4$. \\

Clearly, the product of cardinalities of the blue elements across all chains is an upper bound on $D(n)$. Therefore, 
$$\frac{\log_2(D(n))}{n} <  \frac{233}{720} + \frac{599}{10080} \log_2{3} + \frac{121}{3360} + \sum_{k=5}^{\infty} \frac{\log_2(k)}{2^{k+2}} = 0.4936...$$
$$D(n) < (1.408...)^n$$

\section*{\normalsize 3. Lower Bound for $D(n)$}

For a quick lower bound, consider the LPS given by $D =\{n+1,n+2,\ldots,2n\}$ and observe that for each $q \in (2n/3,n]$, replacing the element $2q \in D$  by $q$ also results in an LPS. Since each replacement is optional, and there are $n/3$ independent decisions, we get a lower bound of $2^{n/3}=(1.259...)^n$. Essentially the same lower bound, attributed to Robert Israel, is mentioned in the OEIS entry of the sequence, referred to in the introduction. \\

We can improve this bound as follows. The idea is to extend the domain of replacement in $D$ from $q \in (2n/3,n]$ to $q \in (n/2,n]$. This cannot be done naively, for at least two reasons. First, it will result in pairs of the form $(t,3t)$ and second, we just proved that $D(n) < (1.408...)^n < 2^{n/2}$. \\

Since all the elements involved are greater than $n/2$, it is clear that the only possible divisor pairs are of the form $(t,3t)$ and $(t,2t)$. We proceed as follows. For each even integer $q \in (n/2,2n/3]$ we consider the quadruple $(q,3q/2,2q,3q)$. These quadruples are made from two $(t,2t)$ pairs, corresponding to $t=q$ and $t=3q/2$. Observe that if $2q$ is replaced by $q$ and $3q$ is left unchanged, we get a $(t,3t)$ pair. Thus we can no longer choose independently whether or not to replace $2q$ and $3q$ with their halves. Therefore, we make both choices simultaneously, and there are three ways to do it, namely $(2q,3q), (q,3q/2)$ and $(2q,3q/2)$. Since there are $n/12$ even integers in $(n/2,2n/3]$, we have $3^{n/12}$ choices. \\

For each integer $q \in (2n/3,n]$ we proceed exactly as before, except when the pair $(q,2q)$ has already occurred as $(3q'/2,3q')$ in our list of quadruples, in which case we discard the pair. Removing $n/12$ such pairs from the original list of $n/3$ pairs, we are left with $2^{n/4}$ choices. Thus, $D(n) > 2^{n/4}3^{n/12} = (1.303...)^n$. 
 
\vskip 10pt

\noindent {\bf Acknowledgement}

\vskip 10pt

I thank Anurag Bishnoi for bringing the problem to my attention.

\end{document}